\documentclass[11pt,letterpaper]{amsart}

\usepackage{amsmath,amssymb,amsfonts,color}
\usepackage[english]{babel}
\usepackage{amsxtra}
\usepackage{mathtools}
\usepackage[normalem]{ulem}
\usepackage[all]{xy}
\usepackage[colorlinks=true, linktocpage, citecolor = blue, linkcolor = blue,backref]{hyperref}
\usepackage{pgf,tikz}
\usetikzlibrary{arrows}

\title[Maximal subgroups in the Cremona group]{Maximal subgroups in the Cremona group}

\newtheorem{theorem}{Theorem}[section]
\newtheorem{lemma}[theorem]{Lemma}
\newtheorem{claim}[theorem]{Claim}
\newtheorem{proposition}[theorem]{Proposition}
\newtheorem{corollary}[theorem]{Corollary}

\newtheorem*{question}{Question}
\newtheorem*{theorem_intro}{Main Theorem}
\theoremstyle{definition}
\newtheorem{definition}[theorem]{Definition}
\newtheorem{construction}[theorem]{Construction}

\theoremstyle{remark}
\newtheorem{remark}[theorem]{Remark}
\numberwithin{equation}{section}

\newcommand{\IA}{\mathbb A}

\newcommand{\IC}{\mathbb C}

\newcommand{\IN}{\mathbb N}
\newcommand{\IP}{\mathbb P}

\newcommand{\IZ}{\mathbb Z}

\newcommand{\cC}{\mathcal{C}}
\newcommand{\cE}{\mathcal{E}}

\newcommand{\cK}{\mathcal{K}}
\newcommand{\cL}{\mathcal{L}}

\newcommand{\cO}{\mathcal{O}}
\newcommand{\cP}{\mathcal{P}}

\newcommand{\cU}{\mathcal{U}}

\newcommand{\Autz}{\mathrm{Aut}^{\circ}}

\newcommand{\Ext}{\mathrm{Ext}}

\DeclareMathOperator{\Hom}{Hom}

\DeclareMathOperator{\Exc}{Exc}

\DeclareMathOperator{\Proj}{Proj}

\DeclareMathOperator{\Aut}{Aut}

\DeclareMathOperator{\Bir}{Bir}
\DeclareMathOperator{\Sym}{Sym}
\DeclareMathOperator{\Chow}{Chow}

\newcommand{\enr}[1]{{\color{teal}#1}} 

\usepackage{tikz}
\usetikzlibrary{cd,arrows,positioning, decorations.pathreplacing}
\tikzset{>=stealth}
\tikzcdset{arrow style=tikz}
\tikzset{link/.style={column sep=1.8cm,row sep=0.23cm}}
\tikzset{link2/.style={column sep=0.4cm,row sep=0.1cm}} 
\tikzset{map/.style={row sep=0em, column sep=0em}}
\tikzset{c/.style={every coordinate/.try}}

\author[A. Fanelli]{Andrea Fanelli}
\address{Andrea Fanelli, 
Univ. Bordeaux, CNRS, Bordeaux INP, IMB, UMR 5251, F-33400 Talence, France.}
\email{andrea.fanelli@math.u-bordeaux.fr}

\author[E. Floris]{Enrica Floris}
\address{Enrica Floris, Universit\'e de Poitiers, Laboratoire de Mathématiques et Applications, UMR 7348 du CNRS, Batiment H3 - Site du Futuroscope, 11 boulevard Marie et Pierre CURIE,
TSA 61125, 86073 Poitiers Cedex 9, France}
\email{enrica.floris@univ-poitiers.fr}

\author[S. Zimmermann]{Susanna Zimmermann}
\address{Susanna Zimmermann, Univ. Paris-Saclay, Institut de mathématiques, CNRS, UMR 8628, F-91405 Orsay, France}
\email{susanna.zimmermann@universite-paris-saclay.fr}

\thanks{
During this project, E.F. and S.Z. were supported by the ANR Project FIBALGA ANR-18-CE40-0003-01. A.F. and E.F. are currently supported by the ANR Project FRACASSO ANR-22-CE40-0009-01. S.Z. was supported the project \'Etoiles Montantes of the R\'egion Pays de la Loire and the Centre Henri Lebesgue and is currently supported by the ERC StG Saphidir and the Institut Universitaire de France.}

\begin{document}
\maketitle

\begin{abstract}
We show that for any $n\geq5$ there exist connected algebraic subgroups in the Cremona group $\Bir(\IP^n)$ that are not contained in any maximal connected algebraic subgroup. Our approach exploits the existence of stably rational, non-rational threefolds.
\end{abstract}

\tableofcontents

\section*{Introduction}

The goal of this work is to elucidate the algebraic structure of higher-dimensional Cremona groups $\Bir(\mathbb{P}^n)$, which are the groups of birational transformations of the $n$-dimensional projective space.
\\It is well-known that $\Bir(\mathbb{P}^1)=\mathrm{PGL}_2(\mathbb{C})$ is an algebraic group, while $\Bir(\mathbb{P}^n)$ with $n\ge 2$ has a much more intriguing group-theoretic nature \cite{CL13,BLZ,LS23,BSY22} and cannot be endowed with the structure of an algebraic group. Thanks to the seminal work by Blanc and Furter \cite{BF13}, we understand the topological obstruction to equip $\Bir(\mathbb{P}^n)$, $n\ge2$, even with the structure of an infinite-dimensional (or ind-)algebraic group. 

In this context, a natural problem consists in studying algebraic groups lying in $\Bir(\mathbb{P}^n)$, $n\ge2$, up to conjugation. 

Demazure formalised in \cite{Dem70} the notion of {rational action} of an algebraic group $G$ on an algebraic variety $X$, i.e.\ of {\it algebraic subgroups of $\Bir(X)$} (see Definition~\ref{def: rational action}). 
After the work of Matsumura \cite{M63}, it is known that $\Bir(X)$ is finite, when $X$ is a variety of general type and, from the view-point of the birational classification of algebraic varieties, we expect $\IP^n$ to lie as far as possible from varieties of general type. It is then natural to interpret the structure of {\it connected} algebraic subgroups of $\Bir(\mathbb{P}^n)$ as a measure of complexity for Cremona groups.

Connected algebraic subgroups of $\Bir(\IP^2)$ have been classified by Enriques \cite{enriques1893sui}: up to conjugacy, they are all contained in (the connected component of the identity $\Aut^\circ$ of) the automorphism group of $\IP^2$ or of (minimal) Hirzebruch surfaces. Moreover the $\Aut^\circ$ of those rational surfaces are all non-conjugate in $\Bir(\IP^2)$. More recently, the classification of maximal finite algebraic subgroups has been completed in \cite{Blanc_alg_subgroups} (see also \cite{DI09,S10}).
 
Maximal connected algebraic subgroups of $\Bir(\IP^3)$ have been classified by Umemura, partially in collaboration with Mukai, in a series of papers \cite{Ume80, Ume82a, Ume82b, MU83, Ume85,Ume88}, see \cite{BFT22, BFT23} for a modern proof using the Minimal Model Program: all connected algebraic subgroups of $\Bir(\IP^3)$ are contained, up to conjugacy, in a maximal one and the full classification of those ones involves several discrete and one continuous families.
The classification of maximal finite algebraic subgroups of $\Bir(\IP^3)$ is not complete, but several results on special classes of finite subgroups of $\Bir(\IP^3)$ have been obtained in the last decade \cite{P11,prokhorov2012simple,BCDP} and it is now clear how modern results in birational geometry can be exploited in the study of Cremona groups \cite{PS16, BLZ} (see also \cite{K20}). 

Any classification in dimension $n\geq4$ is currently unreachable, since we lack fundamental ingredients such as the classification of Fano varieties; partial results have been obtained in \cite{BF20}. 

In this work we are interested in maximal connected algebraic subgroups of the Cremona groups in higher dimensions. In the seminal work \cite{Dem70}, Demazure studied maximal connected algebraic subgroups of $\Bir(\IP^n)$ containing a torus of dimension $n$: his approach originated the study of toric varieties (see also \cite[Section~2.5]{BFT22} for more results on conjugacy classes of tori in $\Bir(\IP^n)$). An interesting feature of connected algebraic subgroup of $\Bir(\IP^2)$ and $\Bir(\IP^3)$ is the following: they are all contained, up to conjugacy, in a maximal one. Blanc asked 10 years ago the following.

\begin{question}
Is every connected algebraic subgroup of $\Bir(\IP^n)$ contained in a maximal one, up to conjugacy?
\end{question}

In \cite{Fong2020, Fong23}, the algebraic subgroups of $\Bir(C\times\IP^1)$, where $C$ is a non-rational curve, are classified up to conjugation. Moreover, Fong shows that if $X$ is a surface of Kodaira dimension $-\infty$, then any algebraic subgroup of $\Bir(X)$ is contained in a maximal algebraic subgroup of $\Bir(X)$ if and only if $X$ is rational. Further results for $\Bir(C\times\IP^n)$, $n\ge2$, have been obtained in \cite{FZ23}.

The main result of this work provides an answer to Blanc's question.

\begin{theorem_intro}\label{thm:main}
For $n\geq 5$, then there exist connected algebraic subgroups of $\Bir(\IP^n)$ which are not contained in any maximal one. 
\end{theorem_intro}

The approach of this work to study this structural question on $\Bir(\IP^n)$ is new and does not depend on any classification, but rather on the nature of rationality in higher dimension.

More concretely, we show the following: let $X$ be the smooth stably rational non-rational threefold of \cite{BCTSSD85}. After birational modification, we may assume that $X$ is endowed with a fibration $c\colon X\to\IP^2$. 
Let $n\geq1$ and consider the projective bundle $\cP_n=\IP_X(\mathcal O_X\oplus c^*\mathcal O_{\IP^2}(n))$ over $X$.
The total space $\cP_n$ has dimension $4$ and since $X\times\IP^m$ is rational for any $m\geq2$ \cite{BCTSSD85, SB02}, the variety $\cP_n\times\IP^m$ is rational for any $m\geq1$. 
We show that for any $n\geq2$ and any $m\geq1$, the connected algebraic subgroup $\Autz(\cP_n\times\IP^m)$ of $\Bir(\IP^{m+4})$ is not contained in any maximal connected algebraic subgroup of $\Bir(\IP^{m+4})$. 

According to the authors' knowledge, it is to date unknown whether $X\times\IP^1$ is rational, i.e.\ whether $\cP_n$ is rational, or not.  Therefore, we are currently unable apply our technique to determine whether  the Main Theorem holds for $n=4$ or not. 

Our construction is inspired by the one used in \cite{FZ23}, where it is shown that for any $n\geq2$ and any curve $C$ of genus $\geq1$, there group $\Bir(C\times\IP^n)$ contains connected algebraic subgroups that are not contained in a maximal connected algebraic subgroup.

\bigskip

\noindent{\bf Acknolwedgements:} We thank Jérémy Blanc, Pascal Fong, Jean-Philippe Furter, Lena Ji, Vladimir Lazic, Andrea Petracci and Sokratis Zikas for interesting discussions. 
We thank János Kollár for pointing out a mistake in a previous version of this work and for his interest in the main problem. He developed in \cite{Kol24a} a new approach which produced new examples and solved the problem in dimension $n=4$.
We also thank Samuel Boissière and Dajano Tossici for the useful discussions involved in fixing the mistake.

\section{Preliminary results}

We work over the field of complex numbers. Varieties are always projective unless stated otherwise.
We refer to \cite{KM98} for the notion of terminality and the basic notions on the minimal model program.

\subsection{Group actions} We recall here some fundamental results on algebraic actions on varieties.
\begin{definition}\label{def: rational action}
Let $Y$ be a variety and let $G$ be an algebraic group.
We say that $G$ acts {\em rationally} on $Y$ if there exists a birational map 
\[
\mu\colon G\times Y\dashrightarrow G\times Y,\ (g,y)\mapsto(g,\mu(g,y))
\] 
that restricts to an isomorphism $U\to V$ on dense open subsets $U,V\subseteq G\times Y$, whose projections onto $G$ are surjective, and such that $\mu(gh,\cdot)=\mu(g,\cdot)\circ\mu(h,\cdot)$ for any $g,h\in G$. If moreover the kernel of the induced homomorphism $G \to \Bir(Y),\ g \mapsto \mu_g$ is trivial, i.e.~if $G$ acts \emph{faithfully} on $Y$, then $G$ is called an {\em algebraic subgroup of $\Bir(X)$}.

The algebraic subgroup $G$ of $\Bir(X)$ is called \emph{maximal} if it is maximal with respect to the inclusion among the algebraic subgroups of $\Bir(X)$.
\end{definition}

Notice that if $W$ is a rational variety and $\psi\colon W\dashrightarrow Y$ a birational map and $G\subseteq\Aut(W)$ an algebraic group, then $G\times Y\to Y$, $(g,y)\mapsto (g,\psi g\psi^{-1}(y))$ is a rational action of $G$ on $Y$ and $G$ is an algebraic subgroup of $\Bir(Y)$; conjugating $G$ by $\psi$ embeds $G$ into $\Bir(Y)$. 

On the other hand, if $G$ is a connected subgroup of $\Bir(Y)$ acting rationally on $Y$,
by the Weil regularisation theorem \cite{Wei55} there is a birational model of $Y$ on which $G$
acts regularly.

\begin{remark}
By \cite{Blanc_alg_subgroups}, any algebraic subgroup of $\Bir(\IP^2)$ is contained in a maximal algebraic subgroup. 
Nevertheless, there are infinite increasing sequences of algebraic subgroups, see \cite[Remark 2.8]{FZ23}. 
\end{remark}

\begin{remark}
It is natural to ask if $\Bir(\IP^n)$ itself can be endowed with a structure of an (ind-)algebraic group. We know this is not possible, thanks to the work \cite{BF13} (see also \cite[Section~2.5]{BFT22} for the construction of the functor $\frak{Bir}_{\IP^n}$).
\end{remark}

\smallskip

We also recall the following two classical facts on regular actions. 
The first follows from \cite[Proposition~2, page~8]{BrionMonastir}.
\begin{lemma}\label{lem:DimMaxopen}
Let $G$ be an algebraic group acting regularly on a projective variety~$X$. Let $n=\max\{\dim(G\cdot x)\mid x\in X\}$ be the maximal dimension of an orbit of $G$. Then, the set $\{x\in X\mid \dim (G\cdot x)<n\}$ is a closed subset of $X$. In particular, the union of orbits of dimension $n$ is a dense open $G$-invariant subset of $X$.
\end{lemma}

The second is the Blanchard's lemma \cite[Proposition~4.2.1]{BSU13}.
\begin{lemma} \label{blanchard}
Let $f\colon X \to Y$ be a proper morphism between varieties such that $f_*(\cO_X)=\cO_Y$. If a connected algebraic group $G$ acts regularly on $X$, then there exists a unique regular action of $G$ on $Y$ such that $f$ is $G$-equivariant.
\end{lemma}

\subsection{Chow varieties}
We refer to \cite{Kollar} for a presentation of Chow varieties, we introduce here the notation and briefly recall some results.\\
Let $X$ be a normal projective variety and $G$ a connected group acting regularly on $X$.
Let ${\rm Chow}(X)$ be the Chow variety of $X$. We recall that ${\rm Chow}(X)$ has countably many irreducible components (cf.\ \cite[Theorem I.3.21(3)]{Kollar}).

If $\mathcal W$ is an irreducible subvariety of ${\rm Chow}(X)$, we denote by $\mathcal U \subseteq \mathcal W \times X$ the universal cycle. Denote by $u\colon \mathcal U \to X$ the natural morphism. If $Z$ is a subvariety of $X$ we denote by $[Z]$ the corresponding point of ${\rm Chow}(X)$. 

Then $G$ acts on every irreducible component of ${\rm Chow}(X)$. Indeed, $G$ preserves every irreducible component as those are countable. If $[Z]$ is a subvariety of $X$ and $g\in G$, then the natural action is given by $g\cdot [Z]=[g(Z)]$. 

\subsection{Ruled varieties}
This section contains some definitions and facts on ruled varieties.
We give first some definitions.

\begin{definition} 
Let $\pi\colon X \to B$ be a morphism between normal projective varieties. One says that $\pi$ is
\begin{enumerate}
\item a $\mathbb{P}^1$-{\it fibration} if its general fibre is a smooth rational curve;
\item a {\it birationally trivial $\mathbb{P}^1$-fibration} if its generic fibre is isomorphic to $\mathbb{P}^1_{\mathbb{C}(B)}$.
\end{enumerate}
Let $\pi\colon X \to B$ be a $\mathbb{P}^1$-fibration. Then one says that $\pi$ is:
\begin{enumerate}
\item[(3)] a {\it standard conic bundle} if $X$ and $B$ are smooth and $\rho(X/B)=1$;
\item[(4)] an {\it embedded conic bundle} if there is a rank 3 vector bundle $\mathcal E$ on $B$ and an embedding $X\hookrightarrow \IP_B(\mathcal E)$ such that
$\pi$ is the restriction of the natural morphism $\rho\colon \IP_B(\mathcal E)\to B$
and $X$ restricted to any fibre of $\rho$ is a conic.
\end{enumerate}
\end{definition}

\begin{remark}\label{rmk:generalities conic fibration}
\
\begin{enumerate}
\item We notice that a birationally-trivial $\mathbb{P}^1$-fibration is a $\mathbb{P}^1$-fibration, and that
$\mathbb{P}^1$-fibrations are also called conic bundles. 

\item Moreover, a fibration is birationally trivial if and only if its general fibre is isomorphic to $\mathbb P^1$ and it admits a birational section.

\item\label{rmk:every stnd cb is embedded} By \cite[Section 1.5]{Sar82}, every  standard conic bundle is embedded. If $G$ is an algebraic group acting regularly the standard conic bundle, then the embedding is equivariant.
\end{enumerate}
\end{remark}

We will often consider {\it projective bundles} of relative dimension $1$, i.e.\ $\mathbb{P}^1$-{\it bundles}, which are projectivisations of a locally free sheaves.
Let $V$ be a projective variety.
If $\mathcal E\to V$ is a rank $r$ vector bundle,
we denote by $\IP(\mathcal E)$ or $\IP_V(\mathcal E)$
the projective bundle of lines in $\mathcal E$
$$\IP_V(\mathcal E)=\Proj(\Sym(\mathcal E))$$
together with the natural morphism $\pi\colon  \IP(\mathcal E)\to V$.

\begin{remark}\label{rem:section of vector bundle}
\
\begin{enumerate}
\item In particular, a surjection $\mathcal E^{\vee}\to\mathcal Q^{\vee}$
determines an embedding $\IP(\mathcal Q)\to \IP(\mathcal E)$ such that $\cO_{\IP(\mathcal Q)}(1)\sim \cO_{\IP(\mathcal E)}(1)\vert_{\IP(\mathcal Q)}$. 
\item By the Noether-Enriques theorem, a smooth $\IP^1$-fibration over a curve is a $\IP^1$-bundle.
\item If $Y$ and $Z$ are smooth, then a $\mathbb{P}^1$-bundle $g\colon Y\to Z$ is a standard conic fibration.
\end{enumerate} 
\end{remark}

\subsection{Rationally connected and non rational threefolds}

We will also need the following statements on rationally connected irrational threefolds.

\begin{proposition}\label{prop:nonratfaut}
 Let $X$ be a rationally connected non-rational threefold. Then $\Autz(X)$ is trivial.
\end{proposition}
\begin{proof}
 Assume by contradiction that $\Autz(X)$ is nontrivial.
 Since $X$ is rationally connected, $\Autz(X)$ is linear and thus contains a 1-parameter subgroup $\Gamma$.
 By \cite[Proposition 2.5.1]{BFT22}, there is an open set $X'$ of $X$ which is of the form $\Gamma\times U$. Since $X$ is rationally connected, any compactification of $U$ is rationally connected.
 Since it is a surface, it is also rational.
 Thus $X$ is birational to $\Gamma\times U$ which is in turn birational to $\IP^1\times\IP^2$, a contradiction.
\end{proof}

\begin{remark}
More generally, one can prove that $\Bir(X)$ contains no connected algebraic subgroups, when $X$ is a rationally connected, non-rational threefold \cite[Corollary~2.5.9]{BFT22}. 
\end{remark}

\begin{remark}
Let $X'$ be a rationally connected threefold. 
Assume $X'$ has a conic bundle structure $X'\to S$. Then there is a birational model
$X$ of $X'$ with a fibration $c\colon X\to \IP^2$ with general fibre $\IP^1$ and sitting in a diagram
\[
\begin{tikzcd}
X\ar[d,"c",swap]\ar[r]&X'\ar[d]\\
\IP^2\ar[r, dashed]&S.
\end{tikzcd}
\]
Indeed, since $X'$ is rationally connected, the surface $S$ is rational.
Let $\IP^2\dasharrow S$ be a birational morphism,
and let $X\to X'\times\IP^2$ be a resolution of the indeterminacies of the induced map $X'\dasharrow \IP^2$. Then the induced morphism $X\to \IP^2$ is the required morphism.

Moreover, if $X'$ is not rational, then the generic fibre of $c$ is not $\IP^1_{\IC(\IP^2)}$.
\end{remark}

\section{From birationally-trivial $\IP^1$-fibrations to $\IP^1$-bundles}
The aim of this subsection is to prove the following statement.

\begin{proposition}\label{pro:proj of rk2 vbdle}
 Let $g\colon Y\to Z$ be a birationally-trivial $\mathbb{P}^1$-fibration. Then there is a 
smooth variety  $\widetilde Z$, a $\IP^1$-bundle $\tilde g\colon \widetilde Y\to\widetilde Z$ and a diagram
 $$
 \xymatrix{
 Y\ar[d]_g\ar@{-->}[r]&\widetilde Y\ar[d]^{\tilde g}\\
 Z\ar@{-->}[r]&\widetilde Z
 }
 $$
 such that all the maps are $\Autz(Y)$-equivariant and the horizontal arrows are birational.
\end{proposition}

\begin{remark}[$G$-equivariant flattening]\label{rem:GRG}
By a result of Gruson and Raynaud \cite{GR71}, if $X,S$ are normal projective varieties and $f\colon X\to S$ is a fibration, there is a birational modification $S'\to S$ such that the base change $f'\colon X\times_S S'\to S'$ is flat.
We notice that, if there is a connected algebraic group $G$ acting on $X$, the construction of $f'$ can be made equivariant for the action of $G$.

Indeed, following \cite[Section 4]{Achinger}, the maximal open set $U$ of $S$
over which $f$ is flat is non-empty and $G$-invariant.
Moreover, $G$ acts on the scheme $\widetilde{\mathcal Q}$ representing the functor $Quot_{\mathcal O_X/X/S}$ by Proposition \ref{prop:App}(1).
Since $G$ is connected, it preserves the (projective) connected components of $\widetilde{\mathcal Q}$.
There is a projective component $\mathcal Q$ of $\widetilde{\mathcal Q}$ such that the sheaf $\mathcal O_X\vert_{f^{-1}U}$ corresponds to a morphism $\psi\colon U\to \mathcal Q$.
The morphism $\psi$ is $G$-equivariant by Proposition \ref{prop:App}(2).
We set $\overline S$ the Zariski closure of the image of $U$ in $\mathcal Q$
and $(p,q)\colon S'\to \overline S\times S$ a $G$-equivariant resolution of the indeterminacy such that $q$ is the composition of blow-ups.
Then the morphism $S'\to \mathcal Q$ corresponds to a flat sheaf on $X\times_S S'$
whose restriction to the preimage of $U$ is $\mathcal O_X\vert_{f^{-1}U}$.
The morphism $S'\to S$ is the required birational modification.
\end{remark}

In \cite[Theorem 1.13]{Sar82}, it is proven that every $\mathbb{P}^1$-fibration is birationally equivalent to a standard conic bundle. 
In our case, result can be made equivariant with respect to the action of a group:

\begin{proposition}\label{prop:standardC}  
 Let  $g\colon Y\to Z$ be a birationally-trivial $\mathbb{P}^1$-fibration. 
Then  there is a standard conic bundle $h\colon V\to S$ and a commutative diagram
$$
\xymatrix{
Y\ar[d]_{g}\ar@{-->}[r]&V\ar[d]^{h}\\
Z&S\ar[l]
}
$$
where the arrows are $\Autz(Y)$-equivariant and $S$ is smooth.
\end{proposition}
\begin{proof}
After a base change and by Remark \ref{rem:GRG}, we can assume that $Z$ is smooth and $g$ is flat.

By the relative Kawamata-Viehweg vanishing theorem \cite[Theorem 3.2.1]{Fujino17}, we have $R^i g_*\cO(-K_Y)=0$ for every $i>0$.
Since $g$ is flat, the rank of $g_*\cO(-K_Y)$ is constant and by \cite[Theorem III 9.9, Corollary III 12.9]{Har77} the sheaf $g_*\cO(-K_Y)$ is locally free of rank three and carries an action of $\Autz(Y)$.

Thus we have a rational map $\sigma\colon Y\dasharrow \IP_Z(g_*\cO(-K_Y))$ over $Z$
which is $\Autz(Y)$-equivariant and an isomorphism onto its image on the open set where $g$ is smooth.
If we set $Y'$ the image of $\sigma$, then $Y'\to Z$ is an embedded conic bundle.
We can now follow the proof of \cite[Theorem 1.13]{Sar82}, by noticing that all the steps are $\Autz(Y)$-equivariant (using Remark~\ref{rmk:generalities conic fibration}\eqref{rmk:every stnd cb is embedded}).
\end{proof}

\begin{proof}[Proof of Proposition \ref{pro:proj of rk2 vbdle}]
By Proposition \ref{prop:standardC}, we may assume that $g$ is a standard conic bundle and $Z$ is smooth.
Then the degeneration divisor $C$ is simple normal crossings.
We prove now that $g$ is a smooth $\IP^1$-fibration, that is, that $C=0$. Let $Z_0$ be a birational section of $g$.
Assume by contradiction that $C$ is non-empty and pick $z\in Supp (C)\setminus Sing (C)$.
Then the fibre over $z$ has two irreducible components $\ell_1$ and $\ell_2$.
Since $Z_0$ is a birational section, we have $Z_0\cdot(\ell_1+\ell_2)=1$.
Since the relative Picard rank is 1, $Z_0$ is relatively ample.
Since $Y$ is smooth, we have  $Z_0\cdot\ell_i\in\IZ$ for $i=1,2$.
A contradiction, because then $Z_0\cdot(\ell_1+\ell_2)\geq2$.

The morphism $g$ is flat because it is equidimensional over a smooth base.
By \cite[Corollary III 12.9]{Har77}  the sheaf $ g_*\mathcal O( Z_0)$ is a rank 2 vector bundle over $Z$.

 Moreover, the natural morphism 
 $Y\to \mathbb P_{Z}(g_*\mathcal O(Z_0))$ is 
 $\Autz(Y)$-equivariant, and an isomorphism on the open set where $g$ is smooth.
 We set therefore $\widetilde Y=\mathbb P_{Z}(g_*\mathcal O(Z_0))$
 and $\tilde g$ the natural morphism.
\end{proof}

\section{Automorphisms of $\IP^1$-bundles}

\subsection{Sections, elementary transformations and automorphism groups}
In this section, we show that invariant sections of projective bundles induce equivariant elementary transformations of $\IP^1$-bundles (Lemma~\ref{lem:CONSTRLEMMA} and Lemma~\ref{lem:invsections}) and we recall the description the automorphism group of projective bundles (Lemma~\ref{lem:auto2}).

\begin{lemma}\label{lem:CONSTRLEMMA}
 Let $V$ be a smooth variety, $\mathcal E\to V$ a rank 2 vector bundle and $\pi\colon\IP(\mathcal E)\to V$ be the induced $\IP^1$-bundle.
 Let $V_0$ be a section of $\pi$ defined by a surjective morphism $\mathcal E^{\vee}\to\mathcal L^{\vee}$. 
 Let $D_1$ be a smooth effective irreducible divisor in $V$. Then the following hold:
 \begin{enumerate}
  \item\label{lem:CONSTRLEMMA1} The sheaf  $\mathcal E_1^{\vee}$ equal to the kernel of the surjection $\mathcal E^{\vee}\to\mathcal L^{\vee}\vert_{D_1}$ is a rank two vector bundle on $V$.
  \item\label{lem:CONSTRLEMMA2} More precisely, if  $\mathcal E^{\vee}$ is an extension
  $$0\to\mathcal M^{\vee}\to\mathcal E^{\vee}\to\mathcal L^{\vee}\to 0$$
  then $\mathcal E_1^{\vee}$ is an extension
  $$0\to\mathcal M^{\vee}\to\mathcal E_1^{\vee}\to\mathcal L^{\vee}(-D_1)\to 0.$$
  \item\label{lem:CONSTRLEMMA3} There is an induced birational map $\psi\colon\IP(\mathcal E)\dasharrow \IP(\mathcal E_1) $
  which factors as $\eta_2\circ\eta_1^{-1}$,
where $\eta_1\colon W\to \IP(\mathcal E)$ is the blow up of the subvariety of $ \IP(\mathcal E)$ defined by $\mathcal E^{\vee}\to\mathcal L^{\vee}\vert_{D_1}$
and $\eta_2\colon W\to \IP(\mathcal E_1)$ is the
contraction of the strict transform of $\pi^{-1}(D_1)$ in $W$.
In particular, $\psi$ is a link 
$$
\xymatrix{
W\ar@{=}[r]\ar[d]_{\eta_1}&W\ar[d]^{\eta_2}\\
\IP(\mathcal E)\ar[d]_{\pi}\ar@{-->}[r]^{\psi}& \IP(\mathcal E_1)\ar[d]^{\pi_1}\\
V\ar@{=}[r]&V
}
$$
 \end{enumerate}

\end{lemma}
\begin{proof}
 We have a diagram of exact sequences
 $$
\xymatrix{
&&0\ar[d]&0\ar[d]&\\
&&\mathcal E^{\vee}_1\ar[d]_{\beta}&\mathcal L^{\vee}(-D_1)\ar[d]&\\
0\ar[r]&\mathcal M^{\vee}\ar[r]^{j}&\mathcal E^{\vee}\ar[r]^{\alpha}\ar[d]_{\alpha_1}&\mathcal L^{\vee}\ar[r]\ar[d]^{\alpha_2}&0\\
&&\mathcal L^{\vee}\vert_{D_1}\ar@{=}[r]&\mathcal L^{\vee}\vert_{D_1}&
}
$$
 Since $\mathcal M^{\vee}=\ker(\alpha)$ and $\mathcal E^{\vee}_1=\ker(\alpha_1\circ\alpha)$, we have an injection $\mathcal M^{\vee}\hookrightarrow\mathcal E^{\vee}_1$.
 Moreover, $\alpha\beta(\mathcal E^{\vee}_1)$ is sent to zero by $\alpha_2$, therefore  $\alpha\beta(\mathcal E^{\vee}_1)$ is contained in $\mathcal L^{\vee}(-D_1)$.
 Via a diagram chase one can prove that the sequence 
  $$0\to\mathcal M^{\vee}\to\mathcal E_1^{\vee}\to\mathcal L^{\vee}(-D_1)\to 0.$$
  is exact.
  Since $\mathcal M^{\vee}$ and $\mathcal L^{\vee}$ have constant rank, the rank of $\mathcal E_1^{\vee}$ is constant as well and we have proved \eqref{lem:CONSTRLEMMA1} and \eqref{lem:CONSTRLEMMA2}.
  
  As for \eqref{lem:CONSTRLEMMA3}, let $U$ be a trivialising set for $\mathcal E^{\vee}$ and $\mathcal E^{\vee}_1$.
  There is a $2\times 2$ matrix $M$ representing the inclusion $\mathcal E^{\vee}_1\vert_U\to \mathcal E^{\vee}\vert_U$.
  If $(e_1,e_2)$ and $(e_1, f_2)$ are local frames for $\mathcal E^{\vee}$ and $\mathcal E^{\vee}_1$ over $U$
  such that $e_1$ is a local frame for $\mathcal M^{\vee}$, then the matrix has the form
  $$M=\left(
  \begin{array}{cc}
   a_{1,1}& a_{1,2}\\
   0& a_{2,2}
  \end{array}
  \right)$$
  where $a_{1,1}\in\Gamma(U,\mathcal O_V^*)$, $a_{2,2}\in\Gamma(U,\mathcal O_V(D_1))$, $a_{1,2}\in\Gamma(U,\mathcal O_V)$.
  The induced map between $\pi^{-1}(U)=U\times\IP^1$ and $\pi_1^{-1}(U)=U\times\IP^1$ is defined by the action of the transposed of $M$.
  Thus we have $$\psi(z,[x_0:x_1])=(z,[a_{1,1}x_0:a_{1,2}x_0+a_{2,2}x_1]).$$
  Without loss of generality we can assume that $a_{1,1}=1$ and multiply by $b=a_{2,2}^{-1}\in\Gamma(U,\mathcal O_V(-D_1))$.
  The section $b$ is a local equation for $D_1$ because $\mathcal E^{\vee}/\mathcal E^{\vee}_1=\mathcal L^{\vee}\vert_{D_1}$ is supported on $D_1$.
  We can assume that there are local analytic coordinates $z=(z_1,\ldots,z_k)$ in $U$ such that $D_1\cap U=\{z_1=0\}$.
  Therefore there are a regular function $f(z)$ on $U$  and a constant $c$ such that $\psi(z,[x_0:x_1])=(z,[c z_1 x_0:z_1 f(z) x_0+x_1]).$
  The indeterminacy locus is thus $D_1\times\{[1:0]\}$. We consider the chart $x_0\neq 0$, set $s=x_1/x_0$ and blow up the ideal $(z_1, s)$.
  The blow up is $$W=\{(z_1,\ldots,z_k,s),[u:v]\vert\;z_1v-su=0\}.$$
  In the chart $u\neq 0$ we have $s=z_1 v/u$.
  Thus on $W$ we extend $\psi$ to a morphism by $\tilde\psi((z_1,\ldots,z_k,s),[u:v])=(z,[c u : f(z)u +v])$. This proves \eqref{lem:CONSTRLEMMA3}.
  \end{proof}

 We give a criterion for the existence of a section fixed by the automorphisms.

\begin{lemma}\label{lem:invsections}
 Let $\mathcal E\to V$ be a rank 2 vector bundle. Assume that $\mathcal E=\mathcal L_1\oplus \mathcal L_2$.
 If $H^0(V,\mathcal L_1\otimes\mathcal L_2^{-1})=\{0\}$,
 then  $\Autz(\mathbb P(\mathcal E))_V$ fixes pointwise the section of $\mathbb P(\mathcal E)\to V$
 corresponding to $\cE^{\vee}\to\mathcal L_2^{-1}$.
\end{lemma}
\begin{proof}
 The projection $\cE^{\vee}\to\mathcal L_2^{-1}$
 induces a section $V_0\to \mathbb P(\mathcal E)$  such that $\mathcal O(1)\vert_{V_0}\sim \mathcal L_2^{-1}$ (see Remark~\ref{rem:section of vector bundle}).
 Viceversa, any section $V_0\to \mathbb P(\mathcal E)$
 with $\mathcal O(1)\vert_{V_0}\sim \mathcal L_2^{-1}$ corresponds to a surjective morphism $\cE^{\vee}\to\mathcal L_2^{-1}$.
 Now, 
 \begin{align*}
 \Hom(\mathcal E^{\vee}, \mathcal L_2^{-1})=&H^0(V, \mathcal L_2^{-1}\otimes \mathcal E) \\
 =&H^0(V, \mathcal L_1\otimes \mathcal L_2^{-1})\oplus H^0(V, \mathcal L_2\otimes \mathcal L_2^{-1})=\mathbb C
 \end{align*}
where the last equality follows from the hypothesis on $\mathcal L_1, \mathcal L_2$.
It follows that $V_0$ is unique with the property that the restriction of $\cO(1)$ to it is $\mathcal L_2^{-1}$ and thus is preserved by the automorphism group.
\end{proof}

We recall the description of the automorphism group of a projective bundle of relative dimension $1$.

\begin{lemma}\label{lem:auto2}
Let $V$ be a smooth variety and $\mathcal E\to V$ a rank-$2$ vector bundle and $\pi\colon\IP(\mathcal E)\to V$ be the induced $\IP^1$-fibration. 
Suppose that $\Aut(\IP(\mathcal E))_V$ fixes a section $V_0$ of $\pi$, given by a surjective morphism $\phi^\vee\colon \cE^\vee\to \cL^\vee$.
Let $\Gamma:=\Gamma(V,\det\cE\otimes \ker(\phi^{\vee})^{\otimes2})$. Then the following hold:
\begin{enumerate}
\item\label{auto2:1} If $\mathcal E$ is decomposable, then $\Aut(\IP(\mathcal E))_V\simeq \Gamma \rtimes \mathbb{G}_m$.
\item\label{auto2:2} If $\mathcal E$ is indecomposable, then $\Aut(\IP(\mathcal E))_V\simeq \Gamma$.
\item\label{auto2:3} If $\Gamma\neq0$, then $V_0$ is the only $\Aut(\IP(\cE))_V$-invariant section. 
\item\label{auto2:4} If $V$ is rationally connected and irrational, then $\Aut(\IP(\cE))\simeq \Aut(\IP(\cE))_V$. 
\item\label{auto2:5} If $V$ is rationally connected and irrational and $\Gamma\neq0$, then the orbits of the $\Autz(\IP(\cE))$-action are included in the fibres of $\pi$ and are either of the form $\pi^{-1}(v)\cap V_0$ for $v\in V$ or the intersection of a fibre $\pi^{-1}(v)$ of $\pi$ with the complement of $V_0$ in $\IP(\cE)$. 
\end{enumerate} 
\end{lemma}
\begin{proof}
We essentially follow \cite[pp.90--92]{Mar71}. 
Let $V=\cup V_i$ be a trivializing cover for $\IP(\cE)$. 
For the morphism $\phi\colon \cL\hookrightarrow\cE$ induced by the surjection surjection $\phi^{\vee}\colon \cE^{\vee}\to\cL^{\vee}$, the image $\mathrm{Im}\phi$ coincides with the annihilator of $\ker(\phi^{\vee})$, which is a hyperplane. By hypothesis, $V_0=\IP(\mathrm{Im}\phi)$ is fixed by $\Autz(\IP(\cE))_V$. If $(v,[x_0:x_1])$ are local coordinates above $V_i$, we can suppose that $V_0$ is given by $x_0=0$. 
Therefore, an automorphism $\varphi\in\Autz(\IP(\cE))_V$ is given by
\[
\varphi_i:=\varphi|_{V_i}=\left(\begin{matrix} \alpha_{i}& s_i\\ 0&1\end{matrix}\right)\in\mathrm{PGL}_2(\cO_V(V_i))
\]
Moreover, the transition functions $\{g_{ij}\}_{i,j}$ of $\IP(\cE)$ are given by
\[
g_{ij}:=\left(\begin{matrix} a_{ij}& c_{ij}\\ 0&1\end{matrix}\right)\in\mathrm{PGL}_2(O_V(V_i\cap V_j))
\]
where the $\{a_{ij}\}_{ij}=\frac{b_{ij}}{d_{ij}}$ and $b_{ij}$ are the transition functions of $\ker(\phi^{\vee})$. 
Notice that the $\varphi_i$ glue to $\varphi\in\Autz(\IP(\cE))_V$ if and only if $g_{ij}\varphi_j=\varphi_ig_{ij}$ for all $i,j$, which is equivalent to 
\[
\alpha_ia_{ij}=a_{ij}\alpha_j,\quad 
a_{ij}s_j+c_{ij}=\alpha_jc_{ij}+s_i \quad\text{for all $i,j$}.
\]
The first condition is equivalent to $\alpha_i=\alpha_j=:\alpha \in \Gamma(V,\cO_V^*)=\mathbb{G}_m$ for all $i,j$.
The second then becomes $s_ja_{ij}-s_i=c_{ij}(\alpha-1)$.

Suppose that $\alpha\neq1$. Then conjugating $g_{ij}$ as follows
\[
\left(\begin{matrix} 1& \frac{s_i}{\alpha-1}\\ 0 &1\end{matrix}\right)
\left(\begin{matrix} a_{ij}& c_{ij}\\ 0&1\end{matrix}\right)
\left(\begin{matrix} 1& -\frac{s_j}{\alpha-1}\\ 0 &1\end{matrix}\right)=\left(\begin{matrix} a_{ij}& 0\\ 0&1\end{matrix}\right)
\]
yields that $\IP(\cE)$ is decomposable, i.e.\ $\IP(\cE)\simeq \IP(\ker(\phi^{\vee})\oplus \cL')$ for some subbundle $\cL'\subset\cE$.
We can then assume that $c_{ij}=0$ and obtain that $s_i=a_{ij}s_j$. Recall that $a_{ij}=\frac{b_{ij}}{d_{ij}}$, where $b_{ij}$ and $d_{ij}$ are respectively the transition functions of $\ker(\phi^{\vee})$ and $\cL'$. 
The $\{d_{ij}^{-1}\}$ define the line bundle $\det\cE\otimes \ker(\phi^{\vee})$, so the $s_i$ glue into a section $s\in\Gamma(V,\det\cE\otimes \ker(\phi^{\vee})^{\otimes2})$. This yields \eqref{auto2:1}. 

If $\alpha=1$, then $s_i=a_{ij}s_j$ and again the $s_i$ glue to a section $s\in\Gamma(V,\det\cE\otimes \ker(\phi^{\vee})^{\otimes2})$ and we obtain \eqref{auto2:2}.

\eqref{auto2:3} If $\Gamma$ is nontrivial, then it is a nontrivial unipotent group it has therefore at most one fixed point on a general fibre of the $\IP^1$-bundle $\pi$. 

\eqref{auto2:4} If $V$ is rationally connected and irrational, then $\Autz(V)$ is trivial by Proposition~\ref{prop:nonratfaut} and hence $\Autz(\IP(\cE))\simeq \Autz(\IP(\cE))_V$. 

\eqref{auto2:5} This follows from \eqref{auto2:3} and \eqref{auto2:4}.
\end{proof}

\begin{remark}\label{rmk:auto2}
Suppose that $V$ admits a conic fibration $c\colon V\to \IP^2$ and that $\cE= \cO_V\oplus c^*\cO_{\IP^2}(n)$. Then $\det\cE\otimes\ker(\phi^{\vee})$ is trivial and in particular, $\Gamma=\Gamma(V,\ker(\phi^{\vee}))\simeq\mathbb{C}[x,y,z]_n$ is the additive group of homogeneous polynomials of degree $n$.
\end{remark}

\subsection{Trivial $\IP^1$-bundles, automorphisms and sections}

The main goal of this section is to prove the following proposition.
\begin{proposition}\label{pro:fix a section}
 Let $g\colon Y\to Z$ be a $\IP^1$-bundle.
 If $\Autz(Y)_Z$ does not fix any section, then $Y=Z\times\IP^1$. 
\end{proposition}

We need two preliminary lemmas.

\begin{lemma}\label{lem:prod}
Let $X,Y$ be smooth projective varieties and $f\colon X\to Y$ be a smooth $\IP^1$-fibration. Assume $f$ has a section $Y_0\subseteq X$.
Then the following are equivalent:
\begin{enumerate}
\item\label{prod:1} $X\cong \mathbb P^1\times Y$ 
\item\label{prod:2} for every general complete intersection curve $\Gamma \subseteq Y$
we have $X_{\Gamma}:= X\times_Y \Gamma\cong  \mathbb P^1\times\Gamma$
and $Y_0\vert_{X_{\Gamma}}$ induces the projection onto $\mathbb P^1$.
\end{enumerate}
\end{lemma}
\begin{proof}
\eqref{prod:1}$\Rightarrow$\eqref{prod:2} is straightforward.
We prove \eqref{prod:2}$\Rightarrow$\eqref{prod:1} by induction on $\dim Y$. If $\dim Y =1$ the claim is true.
Assume thus the claim when the base of the fibration has dimension $n-1$ and assume that $\dim Y=n$.

Let $H\subseteq Y$ be a smooth hyperplane section such that 
\begin{eqnarray}\label{eqn:vanishing}
H^1(Y,f_*\mathcal O_X(Y_0)\otimes\mathcal O_Y(-H))=0
\end{eqnarray}
Let $X_H:=X\times_Y H$. Then \eqref{prod:2} holds in particular for the restriction $f|_{X_H}\colon X_H\to H$. By inductive hypothesis we have $X_H\cong \mathbb P^1\times H$.
By considering the long exact sequence induced by the restriction to $X_H$, we get an exact sequence
$$H^0(X,\mathcal O_X(Y_0))\to H^0(X_H,\mathcal O_{X_H}(Y_0))\to H^1(X,\mathcal O_X(Y_0-X_H))= H^1(X,\mathcal O_X(Y_0-f^*H)).$$
The beginning of the Leray spectral sequence yields an exact sequence

$$0\to H^1(Y,f_*\mathcal O_X(Y_0-f^*H)))\to H^1(X,\mathcal O_X(Y_0-f^*H))\to H^0(Y, R^1 f_*\mathcal O_X(Y_0-f^*H)). $$
By \eqref{eqn:vanishing} we have $H^1(Y,f_*\mathcal O_X(Y_0-f^*H)))=0$.
Moreover the stalk of $R^1 f_*\mathcal O_X(Y_0-f^*H)$ over a closed point has dimension $h^1(\mathbb P^1,\mathcal O(1))=0$.
Therefore $ H^1(X,\mathcal O_X(Y_0-f^*H))=0$ and the restriction map $H^0(X,\mathcal O_X(Y_0))\to H^0(X_H,\mathcal O_{X_H}(Y_0))$ is surjective.

Since $\mathcal O_{X_H}(Y_0)$ is base-point-free by induction hypothesis, the base locus of $\mathcal O_X(Y_0)$ is disjoint from $X_H$.
This holds for every $H$ verifying \eqref{eqn:vanishing}, therefore  $\mathcal O_X(Y_0)$ is base-point-free and defines a morphism $\phi\colon X\to Z$. We want to prove that $Z=\mathbb P^1$.

We now show that $Z$ is a curve. 
Let $y_1,y_2\in Y$ be two distinct points. In order to show that $\dim Z=1$, it suffices to show that $\phi(X_{y_1})=\phi(X_{y_2})$. 
Let $H_1,H_2$ be two hyperplanes in $Y$ satisfying \eqref{eqn:vanishing}. Then $H_1\cap H_2\neq\varnothing$ and we pick $y\in H_1\cap H_2$. For each $i=1,2$, we have a commutative diagram
\[
\begin{tikzcd}
X_{H_i}\ar[r,"\phi|_{X_{H_i}}"]\ar[d,"\iota_{H_i}",hook]& \IP^1\ar[d,"j_{H_i}",hook]\\
X\ar[r,"\phi"] & Z
\end{tikzcd}
\]
where $\iota_{H_i}$ and $j_{H_i}$ are the inclusion. 
For $i=1,2$, we have: since $X_{y_i}, X_y\subset X_{H_i}$ and $\phi|_{X_{H_i}}\colon X_{H_i}\to \IP^1$ is the projection onto the second factor for $i=1,2$, we obtain
\begin{align*}
\phi(X_{y_i})=(\phi\circ \iota_{H_i})(X_{y_i}) 
&= (j_{H_i}\circ\phi|_{X_{H_i}})(X_{y_i})\\
&= (j_{H_i}\circ\phi|_{X_{H_i}})(X_y)=(\phi\circ \iota_{H_i})(X_y)=\phi(X_y).
\end{align*}
Therefore $Z$ is a smooth curve. Since $X_y\cong\IP^1$, we get $\phi(X_y)=Z$, proving that $Z\cong\IP^1$.
\end{proof}

\begin{lemma}\label{cor:fix a section surf}
 Let $S\to C$ be a $\IP^1$-bundle from a projective surface to a curve.
 If $\Autz(S)_C$ does not fix any section, then $S\cong C\times\IP^1$. 
\end{lemma}
\begin{proof}
If $g(C)=0$, then $S$ is isomorphic to a Hirzebruch surface (this is classical, but you may find a mordern proof in \cite[Lemma 2.4.6]{BFT23}). As $\Autz(S)_C$ does not fix any section, we have $S\simeq\IP^1\times\IP^1$. 
If $g(C)\geq 1$, then \cite[Proposition 2.18]{Fong2020} and our assumption imply that either $S\simeq C\times\IP^1$ or that $\min\{\sigma^2\mid \text{$\sigma$ section of $S\to C$}\}>0$. 
In the latter case, $\Autz(S)$ is finite by \cite[Theorem 2(1)]{Mar71} and \cite[Proposition 2.15]{Fong2020}, against our assumption. 
\end{proof}

We are now ready to prove Proposition \ref{pro:fix a section}.

\begin{proof}[Proof of Proposition \ref{pro:fix a section}]
 Let $C$ be a complete intersection curve in $Z$ and $S=Y\times_Z C$ such that the image of the  restriction
 $$\Autz(Y)_Z\to \Autz(S)_C$$
 does not fix any section in $S$.
 Such a curve $C$ exists since $\Autz(Y)_Z$ does not fix any section of $g$. Since $g$ is a $\IP^1$-bundle, the morphism $S\to C$ is a $\IP^1$-bundle.
 By Lemma \ref{cor:fix a section surf} we have $S=C\times\IP^1$.
 By Lemma \ref{lem:prod} we have $Y=Z\times\IP^1$.
\end{proof}

\section{A family of projective bundles over a rationally connected non-rational threefold}
We introduce the projective bundle $\mathcal P_n$, a main player in this article, and prove Proposition~\ref{prop:main} below that characterises $\Autz(\mathcal P_n)$-equivariant birational maps starting from $\mathcal P_n$.

\begin{definition}\label{def:Pn}
Let $X$ be a rationally connected threefold, admitting a fibration $c\colon X\to \IP^2$ with general fibre $\IP^1$. They exist by \cite{BCTSSD85}.
We define $\cP_n=\IP_X(\cO_X\oplus c^*\cO_{\IP^2}(n))$ and $\pi\colon \cP_n\to X$. Let $X_0$ be the section defined by the surjective morphism $\cO_X\oplus c^*\cO_{\IP^2}(-n)\to c^*\cO_{\IP^2}(-n)$.
\end{definition}

Let us recall the properties of $\Autz(\mathcal P_n)$.

\begin{lemma}\label{lem:AutP_n}
Let $X$ be a rationally connected threefold, admitting a fibration $c\colon X\to \IP^2$ with general fibre $\IP^1$.
Assume that $X$ is not rational. Let $\pi\colon\cP_n\to X$ be the $\IP^1$-bundle defined in \eqref{def:Pn}.
The group $\Aut^{\circ}(\cP_n)$ has the following properties:
\begin{enumerate}
\item\label{AutP_n:1} there is an equality $\Aut^{\circ}(\cP_n)=\Aut^{\circ}(\cP_n)_X$;
\item\label{AutP_n:2}  the group $\Aut^{\circ}(\cP_n)$ fixes a section $V_0$ of $\pi$;

\item\label{AutP_n:3} there is an isomorphism $\Aut^{\circ}(\cP_n)\cong \Gamma\rtimes \mathbb{G}_m$, 
where $\Gamma$ is an additive group of dimension $(n+1)(n+2)/2$.
In particular, if $n\geq1$, then $\dim \Autz(\mathcal P_n)\geq 4$;
\item\label{AutP_n:4}  the orbits of the action of $\Autz(\mathcal P_n)$ are included in fibres of $\pi$ and are either of the form $\pi^{-1}(x)\cap X_0$ for $x\in X$ or the intersection of a fibre $\pi^{-1}(x)$ of $\pi$ with the complement of $X_0$ in $\mathcal P_n$.
\end{enumerate}
\end{lemma}
\begin{proof}
The second statement \eqref{AutP_n:2} follows from Lemma \ref{lem:invsections}. 
Statement \eqref{AutP_n:1} follows from the short exact sequence induced by the Blanchard's lemma, or, equivalently, by Lemma \ref{lem:auto2}\eqref{auto2:4} and \eqref{AutP_n:2}.
Statement \eqref{AutP_n:3} follows from Remark~\ref{rmk:auto2}. 
We get \eqref{AutP_n:4} from Lemma \ref{lem:auto2}\eqref{auto2:5}.
\end{proof}

The main goal of this section is to prove the following statement. 

\begin{proposition}\label{prop:main}
Let $n\geq2$ and let $\Phi\colon\mathcal P_n\dasharrow W$ be a birational $\Autz(\mathcal P_n)$-equivariant map.
Then the following hold:
\begin{enumerate}
\item\label{BirStruct:1} $\Autz(\mathcal P_n)$ is a normal subgroup of $\Autz(W)$;
\item\label{BirStruct:2} there are smooth varieties $Y,Z$ and a fibration $Y\to Z$ with generic fibre $\mathbb P^1_{\mathbb C(Z)}$ and a birational 
$\Autz(W)$-equivariant map $\eta\colon W\dasharrow Y$ and a birational map $\varphi\colon X\dashrightarrow Z$ fitting into the following commutative diagram
\[
\begin{tikzcd}
\mathcal P_n\ar[r,dashed,"\Phi"]\ar[d] & W\ar[r,dashed,"\eta"] & Y\ar[d]\\
X\ar[rr,dashed,"\varphi"] && Z
\end{tikzcd}
\]
\end{enumerate}
\end{proposition}

To prove Proposition~\ref{prop:main}, let us fix the following notation and construction. 

\begin{construction}\label{cstr:K}
Notation as in Proposition~\ref{prop:main}.
Since $\Phi$ is $\Autz(\cP_n)$-equivariant, the group $\Autz(\mathcal P_n)$ acts on $W$.
The general orbit has dimension 1, therefore by Lemma \ref{lem:DimMaxopen} the orbits of the action of $\Autz(\mathcal P_n)$ on $W$ have dimension 0 or 1.
Let $\mathcal K_0\subseteq \Chow(W)$ be the subvariety parametrising the orbits of $\Autz(\mathcal P_n)$
and $q_0\colon \mathcal U_0\to \mathcal K_0$ the restriction of the universal family.
Let $\mathcal K\supseteq \mathcal K_0$ be the smallest $\Autz(W)$-invariant closed set in $\Chow(W)$
and $u\colon \mathcal U\to \mathcal K$ the restriction of the universal family. 
\[
\begin{tikzcd}
\mathcal U_0\ar[r, hook]\ar[d,"u_0",swap]&\mathcal U\ar[d,"u"]\ar[r,"e"]& W\\
\mathcal K_0\ar[r, hook]&\mathcal K&
\end{tikzcd}
\]
There are evaluation morphisms $e\colon \mathcal U\to W$ and $e_0=e\vert_{\mathcal U_0}\colon \mathcal U_0\to W$. Those are $\Autz(W)$-equivariant and $\Autz(\mathcal P_n)$-equivariant respectively.
We notice that $e_0$ is birational by Lemma~\ref{lem:AutP_n}. 
\end{construction}

In the following lemmas we follow the notation from Construction \ref{cstr:K} and Proposition \ref{prop:main}.

\begin{lemma}\label{lem:norm}
The group $\Autz(\cP_n)$ is a normal subgroup of $\Autz(W)$ if and only if $\mathcal K_0=\mathcal K$. 
\end{lemma}
\begin{proof}
The group $\Autz(\cP_n)$ is a normal subgroup of $\Autz(W)$ if and only if $\Autz(W)$ permutes the $\Autz(\cP_n)$-orbits in $W$. This is the case if and only if $\Autz(W)$ preserves $\mathcal K_0$. By minimality of $\mathcal K$, this is equivalent to $\mathcal K_0=\mathcal K$. 
\end{proof}

\begin{lemma}\label{lem:contr}
If $\Phi$ does not contract $X_0$, then $\mathcal K_0=\mathcal K$. 
\end{lemma}
\begin{proof}
By Lemma~\ref{lem:AutP_n}, the only $\Autz(\cP_n)$-invariant proper closed subvarieties of $\cP_n$ are union of fibres of $\pi$ or $X_0$. If $X_0$ is not contracted by $\Phi$, then there exists an open nonempty subset $U\subset X$ such that $\Phi|_{\pi^{-1}(U)}$ is an isomorphism. 

Suppose that $\mathcal K_0\neq \mathcal K$ and let $M$ be the pull-back by $\pi\colon \cP_n\to X$ of an ample divisor on $X$.  
Let $[\Gamma]\in\mathcal K\setminus \mathcal K_0$ and $[\Gamma_0]\in\mathcal K_0
$ be classes of curves $\Gamma,\Gamma_0$ on W such that $\Gamma_0\subset\Phi(\pi^{-1}(U))$ and such that $\Gamma$ is not in the exceptional locus of $\Phi^{-1}$. 

We denote by $(p,q)\colon\hat W\to \cP_n\times W$ an $\Autz(\cP_n)$-equivariant resolution of $\Phi$. By Lemma~\ref{lem:AutP_n}, the conic bundle $\pi\colon \cP_n\to X$ has a unique $\Autz(\mathcal P_n)$-invariant section $X_0\subset\cP_n$, and by $\hat{X}_0$ we denote its strict transform in $\widehat W$.

Let $\widehat\Gamma_0$ be the strict transform of $\Gamma_0$ in $\widehat{W}$ and let $\widehat\Gamma$ be an irreducible curve in $\hat W$ such that $q_*(\widehat\Gamma)=\Gamma$. Notice that $\widehat\Gamma$ is not contracted by $\pi\circ p$, because $[\Gamma]\notin\mathcal K_0$. 
Then

\begin{claim}\label{cla:pback of curves}
\[
p^*M\cdot\widehat\Gamma_0=p^*M\cdot\widehat\Gamma+p^*M\cdot C
\]
for some curve $C$ in $\hat W$. 
\end{claim}
Assuming the claim, we finish the proof.
The left-hand side is zero, because $[\Gamma_0]\in\mathcal K_0$, while $p^*M\cdot\widehat\Gamma>0$ and $p^*M\cdot C\geq0$. This is impossible, so it follows that $\mathcal K_0=\mathcal K$. 

We are left with the proof of Claim \ref{cla:pback of curves}.
Let $(a,b)\colon\widehat{\cU}\to \cU\times\widehat W$ be a resolution of the indeterminacies of $\cU\dasharrow\widehat W$ and let $\hat u\colon \widehat{\cU}\to \cK$ be the induced fibration. 
Let $\cC$ be an irreducible curve in $\cK$ such that $[\Gamma],[\Gamma_0]\in \cC$.
Let $S$ be the component of dimension 2 of $\hat u^{-1}(\cC)$ surjecting onto $\cC$.
Then $b_*\hat u^*[\Gamma_0]=\widehat \Gamma_0$ and there is an effective curve $C$ such that $b_*\hat u^*[\Gamma]=\widehat \Gamma +C$. The claim follows as  $b_*\hat u^*[\Gamma_0]\equiv b_*\hat u^*[\Gamma]$.
\end{proof}

\begin{lemma}\label{lem:invariance of ideal}
Suppose that $\Autz(\cP_n)$ is a not normal subgroup of $\Autz(W)$ and that 
every $\Autz(W)$-equivariant desingularisation of $W$ extracts $X_0$.
Let $[\Gamma_0]\in \mathcal K_0$ be a general point, 
let $G\subset\Autz(W)$ be a $1$-parameter subgroup and $g\in G$ a general element. 
Let $\widetilde{g\Gamma_0}$ be the strict transform of $g\Gamma_0$ in $\mathcal P_n$.
Then $\widetilde{g\Gamma_0}\cap  X_0$ is a non-empty finite set. 
\end{lemma}
\begin{proof}
Let $\widehat W\to W$ be an $\Autz(W)$-equivariant desingularisation. Then $\widehat W\to W$ extracts $X_0$.
We denote by $\hat X_0$ the strict transform of $X_0$ in $\widehat W$. Then the induced birational map $\mathcal P_n\dashrightarrow \widehat W$ is an isomorphism at the generic point of $X_0$ and induces a birational map $X_0\dashrightarrow \widehat X_0$.
Let $[\Gamma_0]\in\mathcal K_0$ be the class of a general curve $\Gamma_0$ such that its strict transform $\widehat\Gamma_0$ in $\widehat W$ meets $\widehat X_0$ in a point lying in the open set where $\widehat W\dasharrow \mathcal P_n$ is an isomorphism.
Then for general $g\in G$, the curve $g\widehat\Gamma_0$  meets $\widehat X_0$ in a point lying in the open set where $\widehat W\dasharrow \mathcal P_n$ is an isomorphism.
Let $\widetilde{g\Gamma_0}$ be the strict transform of $g\Gamma_0$ in $\mathcal P_n$.
Then $\widetilde{g\Gamma_0}\cap  X_0$ is non-empty. 
\end{proof}

\begin{lemma}\label{lem:not invariance of ideal}
Suppose that $\Autz(\cP_n)$ is a not normal subgroup of $\Autz(W)$.
Then there is an $\Autz(W)$-equivariant desingularisation of $W$ which does not extract $X_0$.

\end{lemma}
\begin{proof}
We prove the statement by contradiction. 
Suppose that all $\Autz(W)$-equiva-\\riant desingularisation of $W$ extract $X_0$.\\ 
Since $\Autz(\cP_n)$ is  not normal in $\Autz(W)$, by Lemma \ref{lem:norm} we have $\mathcal K_0\subsetneq\mathcal K$. Then there is a $1$-parameter subgroup $G\subset\Autz(W)$
with $G\not\subseteq\Autz(\cP_n) $ and such that for a general $g\in G$, for a general $[\Gamma_0]\in\cK_0$, we have $[g\Gamma_0]\not\in\cK_0$. 
Let $C$ be the strict transform of $g\Gamma_0$ in $\mathcal P_n$. 
Lemma~\ref{lem:invariance of ideal} 
implies that $C\cap  X_0$ is non-empty. 
Let $H\subset\Autz(\cP_n)$ be an additive $1$-parameter subgroup, set $C_t:=tC$, $t\in H$, and consider the pencil $\{C_t\}_{t\in H}$. 
The pencil defines a morphism $\mu\colon \IA^1\times\IP^1\to \cP_n$. 
Let $F$ be the normalisation of the Zariski-closure of $\mu(\IA^1\times\IP^1)$ and $n\colon F\to \mathcal P_n$ the induced morphism. The image $\pi(C)$ is a curve because $[g\Gamma_0]\not\in\cK_0$. Let $n\colon D\to \pi(C)$ be the normalisation of $\pi(C)$. By abuse of notation, the strict transform of $C_t$ (resp. $C$) on $
F$ is denoted by $C_t$ (resp. $C$) as well, as no confusion will arise. 

Notice that $F$ is smooth and that it is a $\IP^1$-bundle above $D$. 
Let $\theta\colon S\to F$ be a minimal resolution of the base-locus of the pencil $\{C_t\}_{t\in H}$. Then there is a conic fibration $u\colon S\to \IP^1$, whose fibres are the strict transforms $\overline{C_t}$ of $C_t$, such that the following diagram commutes.
\[
\begin{tikzcd}
\IA^1\times\IP^1\ar[d,"u"]\ar[r,hook]\ar[rrr, bend left=20,"\mu"]&S\ar[r,"\theta"]\ar[d,"u"]&F\ar[r,"n"]\ar[d,"\pi"]&\cP_n\ar[d,"\pi"]\\
\IA^1\ar[r,hook]&\IP^1&D\ar[r,"n"]&X
\end{tikzcd}
\]
Notice that since $\Gamma_0$ is rational, so is $C=g\Gamma_0$ and thus $D\simeq\IP^1$. 
By abuse of notation, the strict transform of $X_0$ on $F$ will be denoted by $X_0$ as well, as no confusion will arise. Notice that $H$ fixes $X_0$ pointwise as $H\subset\Autz(\cP_n)$, so all the points in $C\cap X_0$ are base-points of the pencil $\{C_t\}_{t\in H}$. 

The $H$-action on $F$ lifts to $S$ and permutes non-trivially the fibres of $u$, since any fibre of $\pi$ intersects $C$ in only finitely many points.  
By the Blanchard's lemma (\cite{Bla56}, \cite[Proposition~4.2.1]{BSU13}), $H$ acts on $\IP^1$ (base of the fibration $u$) and, since it is additive, it fixes exactly one point  $[\infty]\in\IP^1$.
Moreover, a curve $E$ contained in the exceptional locus of $\theta$ is either contained in a fibre of $u$ or it is a section of $u$, and the latter is the case if and only if $E$ is the exceptional divisor of a point that is blown up last, or, equivalently, a $(-1)$-curve.
\medskip

{\bf Claim.} $C_t$ is a section of $\pi\colon F\to D$ if $t\in H$ is general.\\
{\em Proof.} Let $f$ be a fibre of $\pi$, which is disjoint from the exception locus of $\theta$ and such that $H$ acts non-trivially on $f$. The fibre being disjoint from the exceptional locus implies that the pullback $\theta^*f$ and the strict transform $\bar f$ of $f$ in $S$ coincide.
The action of $H$ being non-trivial implies that $\bar f$ is a section of $u$.
Indeed, the restriction $u\colon \bar f\to \IP^1$ is surjective and $H$-equivariant, therefore the ramification and branch locus are preserved by $H$.
But those are both supported on at most one point and by Hurwitz formula this is possible only if $u$ has degree 1. 
Let $t$ be such that $\overline C_t$ is irreducible.
Thus we get $$f\cdot C_t=\theta^*f\cdot\overline C_t=\bar f\cdot\overline C_t=1$$
where the first equality is the projection formula, the second is because the pullback $\theta^*f$ coincides with the strict transform $\bar f$, and the third because $\bar f$ is a section of $u$. This finishes the proof that $C_t$ is a section of $\pi\colon F\to D$ if $t\in H$ is general.
\smallskip

In particular, $C=C_0$ is a section of $\pi$. 
It follows that every fibre of $\pi$ meets the base-locus in $F$ in at most one point.
Since $C\cap X_0$ is non empty, we pick a point $y_0\in C\cap X_0$, we set $0=\pi(y_0)\in D$ and $f_0=\pi^{-1}(0)$.
Let $f_0=f_{t_0},f_{t_1},\dots,f_{t_l}$ be the fibres of $\pi$ meeting the base-locus of the pencil $\{C_t\}_{t\in H}$, and let $\overline{f_{t_i}}$ (resp. $\overline{f_0}$) be the strict transform of $f_{t_i}$ (resp. $f_0$) in $S$. Since $C_t$ is a section of $\pi$, 
the $f_{t_i}$ are not contained in the union $\cup_{t\in H}C_t=\theta(\IA^1\times\IP^1)$, and more precisely the intersection of $f_{t_i}$ with $\cup_{t\in H}C_t$ coincides with one point. Therefore each $f_{t_i}$ is contained in $\theta(u^{-1}[\infty])$. 

For each $i\geq0$, denote by $E_{ij}$ the irreducible components of $\theta^{-1}(y_i)$ that are not sections of $u$ and by $E_i^{sec}$ the unique irreducible component that is a section of $u$. 
Since $\Exc(\theta)$ is preserved by $H$, we have $\theta(E_{ij})=[\infty]$ for every $i,j$.
Then 
\[
\theta^*C_0=\overline C_0+\sum a_{ij} E_{ij}+\sum_{i=0}^l a_i E_i^{sec}
\]
for some integers $a_{ij},a_i\geq0$. Denote by $\overline C_{\infty}$ the fibre of $u$ above $[\infty]$. 
Since all $E_{ij}$ are contained in $\overline C_{\infty}$ and the $E_i^{sec}$ are sections of $u$, we have 
\[
\overline C_{\infty}\cdot\theta^*C_0=C_t\cdot\theta^*C_0=\sum_{i=0}^l a_i,\quad\text{for any $t\in H$}.
\] 
We also have
\[
\overline C_{\infty}=\theta^*(\alpha X_0+\sum_{i=0}^l \beta_i f_{t_i})-\sum_{i=0}^lb_i E_i^{sec}+C'
\]
for some integers $\alpha,\beta_i,b_i\geq0$ and some effective divisor $C'$ having no common component with $\sum E_i^{sec}$. In fact, $b_i\geq\beta_i\geq1$ for all $i\geq1$, because $E_i^{sec}\subseteq\mathrm{supp}(\theta^*f_{t_i})$ and $\bar f_{t_i}$ is contained in $\overline{C_{\infty}}$ for $i\geq1$. 
Furthermore,  $b_0\geq\beta_0+\alpha\geq2$, because $E_0^{sec}\subset\mathrm{supp}(\theta^*f_{0})\cap \mathrm{supp}(\theta^*X_0)$ and $\bar f_0$ is contained in $\overline{C_{\infty}}$. We compute
\begin{align*}
\sum_{i=0}^l a_i
&=\bar C_{\infty}\cdot\theta^*C_0=\bar C_{\infty}\cdot\left(\sum_{i=0}^l a_i E_i^{sec}\right)\\
&=\left(\theta^*(\alpha X_0+\sum_{i=0}^l  \beta_i f_{t_i})-\sum_{i=0}^l b_i E_i^{sec}+C'\right)\left(\sum_{i=0}^l a_i E_i^{sec}\right)\\
&\stackrel{(E_i^{sec})^2=-1}=\sum_{i=1}^l b_ia_i+C'\cdot \left(\sum_{i=1}^l a_i E_i^{sec}\right)\geq \sum_{i=0}^l b_ia_i\geq a_0+\sum_{i=0}^l a_i.
\end{align*}
where the last inequality holds because $b_0\geq2$ and $b_i\geq1$ for $i\geq1$.
It follows that $a_0=0$, which contradicts $f_0$ containing a base-point of $\{C_t\}_{t\in H}$. 
This proves that there is an $\Autz(W)$-equivariant desingularisation of $W$ which does not extract $X_0$.
\end{proof}

\begin{lemma}\label{lem:contr0}
Suppose that $\Autz(\cP_n)$ is a not normal subgroup of $\Autz(W)$.
Then $\Phi$ does not contract $X_0$. 
\end{lemma}
\begin{proof}
Suppose that $\Phi$ contracts $X_0$. 
By Lemma~\ref{lem:not invariance of ideal}, there exists $\mu\colon \widetilde{W}\to W$ an $\Autz(W)$-equivariant desingularisation of $W$ which does not extract $X_0$. 

We denote by $(p,q)\colon\widehat W\to \cP_n\times W$ an $\Autz(\cP_n)$-equivariant resolution of the indeterminacy of $\Phi$. By Lemma~\ref{lem:AutP_n}\eqref{AutP_n:3}, the conic bundle $\pi\colon \cP_n\to X$ has a unique $\Autz(\mathcal P_n)$-invariant section $X_0\subset\cP_n$, and by $\widehat{X}_0$ we denote its strict transform in $\widehat W$. 

 Let $(\bar{p},\bar{q})\colon\overline W\to \widehat W\times\widetilde W$ be a resolution of the indeterminacies of $\mu^{-1}q\colon \widehat{W}\dashrightarrow\widetilde W$, such that $\bar q$ is a composition of blow-ups of smooth centres. 
\[
\begin{tikzcd}
&&\overline{W}\ar[dl,"\bar{p}",swap]\ar[dr, "\bar{q}"]&\\
&\widehat{W}\ar[dl,"p",swap]\ar[dr,"q"]&&\widetilde W\ar[dl,"\mu",swap]\\
\cP_n\ar[rr,"\Phi",dashed]&& W&
\end{tikzcd}
\]
Then the strict transform $\overline{X}_0$ in $\overline W$ is among the exceptional divisors of those blow-ups. 
It follows that $X_0$ is birational to $\IP^k\times Z$ with $k\in\{1,2,3\}$ and $\dim Z=3-k$. Since $\dim Z\leq 2$ and $Z$ is rationally connected, $Z$ is rational and so is $X_0$. This contradicts the hypothesis that $X$ is not rational. 
\end{proof}

\begin{proof}[Proof of Proposition~\ref{prop:main}]
We first prove that $\mathcal K_0=\mathcal K$.
If not, by Lemma \ref{lem:norm}, the group $\Autz(\mathcal P_n)$ is not normal in $\Autz(W)$.
By Lemma \ref{lem:contr0} the map $\Phi$ does not contract $X_0$.
This is a contradiction with Lemma \ref{lem:contr}.

Therefore $\mathcal K_0=\mathcal K$ and $\Autz(\mathcal P_n)$ is normal in $\Autz(W)$ by Lemma \ref{lem:norm}.
Let $\overline{\mathcal U}$ and $\overline{\mathcal K}$ be $\Autz(W)$-equivariant compactifications of $\mathcal U$ and $\mathcal K$
such that there is a morphism $u\colon \overline{\mathcal U}\to \overline{\mathcal K}$ extending $u$.
The map $e\colon \overline{\mathcal U}\dasharrow W$ is birational and $\Autz(W)$-equivariant, see Construction~\ref{cstr:K}.
We set $Y=\overline{\mathcal U}$ and $Z=\overline{\mathcal K}$.

Moreover, $Z$ is birational to $X$ because $\mathcal K_0=\mathcal K$ parametrises the 1-dimensional orbits of $\Autz(\mathcal P_n)$.
\end{proof}

\section{$\Autz(\mathcal P_n)$ is not contained in a maximal subgroup of $\Bir(\mathcal P_n)$}
The aim of this section is to show in Theorem~\ref{thm:main1new} that $\Autz(\mathcal P_n)$ is not contained in a maximal connected algebraic subgroup of $\Bir(\mathcal P_n)$ if $n\geq2$.

\begin{proposition}\label{pro:linkIInew}
Let $V$ be a smooth variety of dimension at least 3. Let $\mathcal E\to V$ be a rank 2 vector bundle.
Assume that $\Autz(\IP(\mathcal E))_V$ contains a non-trivial  additive group and fixes a section of $\pi\colon \IP(\mathcal E)\to V$. 
Then there is a rank 2 vector bundle $\mathcal E_1\to V$ and an $\Autz(\IP(\cE))_V$-equivariant birational map $\IP(\cE)\dashrightarrow \IP(\cE_1)$ over $V$
such that $\Autz(\IP(\mathcal E))_V\subsetneq \Autz(\IP(\mathcal E_1))_V$.
\end{proposition}

\begin{proof}
  From the two hypotheses on $\Autz(\IP(\mathcal E))_V$, it follows that there is a unique section $V_0$
 fixed by $\Autz(\IP(\mathcal E))_V$.
 This section corresponds to the data of a line bundle $\mathcal L$ on $V$ and a surjective morphism $\mathcal E^{\vee}\to\mathcal L^{\vee}$ (see Remark~\ref{rem:section of vector bundle}).
 Let $\mathcal M^{\vee}$ be the kernel of $\mathcal E^{\vee}\to\mathcal L^{\vee}$.
 Then $\mathcal M^{\vee}$ is a rank 1 torsion-free sheaf on $V$ and it is locally free by \cite[Proposition 1.9]{Har80}.

We have $\Autz(\IP(\mathcal E))_V\simeq \Gamma\rtimes G$,  for some non-trivial additive group $\Gamma$ and $G=\mathbb G_m$ or $G=\{1\}$, see Lemma \ref{lem:auto2}.
  Let $\overline D$ be a very ample divisor on $V$ such that 
 \begin{itemize}
  \item there is a smooth element $D_1\in |\overline D|$;
  \item the line bundle $\mathcal M^{\vee}\otimes \mathcal L(\overline D)$ is very ample, so that 
  $\Ext^1_V(\mathcal L^{\vee}(-\overline D),\mathcal M^{\vee})=H^1(\mathcal M^{\vee}\otimes \mathcal L( \overline D))=0$; and 
  \item $\dim \Gamma<\dim H^0 (V,\mathcal M^{\vee}\otimes \mathcal L( \overline D)).$
 \end{itemize}
 Consider now the kernel $\mathcal E_1^{\vee}$ of the surjection $\mathcal E^{\vee}\to\mathcal L^{\vee}\vert_{D_1}$.
By  Lemma~\ref{lem:CONSTRLEMMA}\eqref{lem:CONSTRLEMMA2}, the sheaf $\mathcal E_1^{\vee}$ is an extension of $\mathcal L^{\vee}(-\overline D)$ and $\mathcal M^{\vee}$ and since $\Ext^1_V(\mathcal L^{\vee}(-\overline D),\mathcal M^{\vee})=0$
 we have $\mathcal E_1^{\vee}\cong \mathcal M^{\vee}\oplus \mathcal L^{\vee}(-\overline D)$.
 Therefore $\mathcal E_1$ is decomposable. 
By Lemma~\ref{lem:invsections}, the group $\Autz(\IP(\mathcal E_1))_V$ fixes the section corresponding to
 $\mathcal E_1^{\vee}\to \mathcal L^{\vee}(-\overline D)$. Then by Lemma \ref{lem:auto2}\eqref{auto2:1}
  $\Autz(\IP(\mathcal E_1))_V\simeq H^0 (V,\mathcal M^{\vee}\otimes \mathcal L( \overline D)) \rtimes \mathbb{G}_m$
 
 Moreover, by the invariance of $V_0$, the link or birational map $\psi\colon \IP(\cE)\dashrightarrow\IP(\cE_1)$ obtained by Lemma~\ref{lem:CONSTRLEMMA}\eqref{lem:CONSTRLEMMA3}
 is $\Autz(\IP(\mathcal E))_V$-equivariant. 
Since $\Gamma\subsetneq H^0 (V,\mathcal M^{\vee}\otimes \mathcal L( \overline D))$ by assumption on $\overline D$, we have $\Autz(\IP(\mathcal E))_V\subsetneq \Autz(\IP(\mathcal E_1))_V$.
 This proves the claim.
\end{proof}

We are now ready to prove the main theorem of this section.

\begin{theorem}\label{thm:main1new}
Let $n\geq 2$ be a positive integer.
Let $X$ be a non-rational and rationally connected variety carrying a non-trivial conic bundle structure and admitting a $\IP^1$-fibration $c\colon X\to \IP^2$.
Set $\mathcal P_n=\IP_X(\mathcal O_X\oplus c^*\mathcal O_{\IP^2}(n))$.
The group $\Autz(\mathcal P_n)$ is not contained in a maximal group of $\Bir(\mathcal P_n)$.
More precisely, for every variety $W$, for every $\Autz(\mathcal P_n)$-equivariant birational map $W\dasharrow \mathcal P_n$, there is a variety $Y$ and an  $\Autz(\mathcal P_n)$-equivariant birational map
$W\dasharrow Y$ with $\Autz(W)\subsetneq\Autz(Y)$.
\end{theorem}
\begin{proof}
Assume that $\Autz(\mathcal P_n)$ is contained in a connected algebraic subgroup $H$ of $\Bir(\mathcal P_n)$ acting rationally on $\mathcal P_n$.
We will prove that there is a connected algebraic subgroup $G$  of $\Bir(\mathcal P_n)$ acting rationally on $\mathcal P_n$ such that $H\subsetneq G$.

By the Weil regularisation theorem \cite{Wei55}, there is a variety $W$ birational to $\mathcal P_n$
such that $H\subseteq \Autz(W)$.
By Proposition \ref{prop:main}, 
there are smooth varieties $Y,Z$ and a fibration $g\colon Y\to Z$ with generic fibre $\mathbb P^1_{\mathbb C(Z)}$, a birational 
$\Autz(W)$-equivariant map $W\dasharrow Y$ and a birational map $X\dashrightarrow Z$ fitting into a commutative diagram
\[
\begin{tikzcd}
\mathcal P_n\ar[r,dashed]\ar[d] & W\ar[r,dashed] & Y\ar[d,"g"]\\
X\ar[rr,dashed] && Z
\end{tikzcd}
\]
By Proposition \ref{pro:proj of rk2 vbdle}, we can assume that $Z$ is smooth and $g$ is a $\IP^1$-bundle.

We claim the following.
\begin{claim}\label{cla:sectionfixed}
 The $\IP^1$-bundle $g$ has an $\Autz(Y)$-equivariant section.
\end{claim}

Assuming the claim, we finish the proof. Since $Z$ is birational to $X$, by Proposition \ref{prop:nonratfaut} we have $\Autz(Y)=\Autz(Y)_Z$.
Since $\Autz(\cP_n)\subseteq\Autz(Y)$, there is a non-trivial additive subgroup of $\Autz(Y)=\Autz(Y)_Z$.
By the Claim \ref{cla:sectionfixed}, the group $\Autz(Y)$ fixes a section of $g$. Therefore, by Proposition \ref{pro:linkIInew}, there is a $\IP^1$-bundle $\IP(\mathcal E_1)\to Z$
such that $\Autz(Y)=\Autz(Y)_Z\subsetneq \Autz(\IP(\mathcal E_1))_Z=\Autz(\IP(\mathcal E_1))$. 
We may set $G=\Autz(\IP(\mathcal E_1))$.

\smallskip

We are left with the proof of the Claim \ref{cla:sectionfixed}: by Proposition \ref{pro:fix a section}, if there is no $\Autz(Y)$-equivariant section then $Y=Z\times\IP^1$.
Then we would have $\Autz(Y)\cong \mathrm{PGL}_2(\mathbb{C})$. But this contradicts the fact that $\Autz(Y)$ contains $\Autz(\cP_n)$,
which has dimension at least 4 by Lemma~\ref{lem:AutP_n}\eqref{AutP_n:3}.
\end{proof}

\section{Proof of Main Theorem}

We start with some preliminary lemmas on birational map from products with special properties.

\begin{lemma}\label{lem:nef of product is product of nef}
Let $X_1$ and $X_2$ be normal projective varieties such that $h^1(X_i,\mathcal O_{X_i})=0$ for $i=1,2$ and let $p_i\colon X_1\times X_2\to X_i$ be the projection onto $X_i$. Then 
\[
\mathrm{Nef}(X_1\times X_2)=p_1^*\mathrm{Nef}(X_1)\oplus p_2^*\mathrm{Nef}(X_2).
\] 
\end{lemma}
\begin{proof}
By \cite[Exercise III 12.6(b)]{Har77} and since $h^1(X_i,\mathcal O_{X_i})=0$ for $i=1,2$, we have $\mathrm{Pic}(X_1\times X_2)=p_1^*\mathrm{Pic}(X_1)\oplus p_2^*\mathrm{Pic}(X_2)$. 
Let $L\subset\mathrm{Nef}(X_1\times X_2)$ and write $L=p_1^*D_1+p_2^*D_2$ for some $D_i\in\mathrm{Pic}(X_i)$, $i=1,2$. 
Let $C_1\subset X_1$ be a curve and consider the curve $\hat{C}=C_1\times\{x_2\}\subset X_1\times X_2$. 
Then $0\leq L\cdot \hat C=D_1\cdot C$ and hence $D_1$ is nef. The same argument shows that $D_2$ is nef. 
\end{proof}

\begin{proposition}\label{prop:bir maps between products}
Let $P$ be a smooth projective variety and $Y$ a homogeneous variety with $\rho(Y)=1$ and let $\varphi\colon P\times Y\dashrightarrow Q$ be an $\Autz(P\times Y)$-equivariant birational map. Then $Q\simeq P'\times Y$, where $P'$ is projective and $\varphi=(\varphi_1,\varphi_2)$ with $\varphi_1\colon P\dashrightarrow P'$ birational and $\varphi_2\colon Y\to Y$ an isomorphism. 
\end{proposition}
\begin{proof}
Let $(p,q)\colon W\to P\times Y\times Q$ be a functorial resolution of the indeterminacies of $\varphi$ such that $p$ is a composition of blow-ups of smooth centres. 
Since $p$ is $\Autz(P\times Y)$-equivariant and $Y$ is homogeneous, the morphism  $p$ blows up centres that are products of the form $C_i\times Y$. 
It follows that $W\simeq\hat P\times Y$ and that $p=(p_{\hat P},p_Y)$, with $p_{\hat P}$ birational and $p_Y$ an isomorphism. 
 
Since $q$ is a birational morphism, it is induced by a Cartier divisor $D$ on $\hat P\times Y$ that is big and nef. 
By Lemma~\ref{lem:nef of product is product of nef}, we can write $D=p_{\hat P}^*D_1+p_Y^*D_2$ and $D_1\in\mathrm{Nef}(\hat P)$ and $D_2\in\mathrm{Nef}(Y)$. 
Then $q$ is of the form $q=(f_1,f_2)\colon \hat P\times Y\to P'\times Y'$, where 
$f_1\colon\hat P\to P'$ is defined by $D_1$ and $f_2\colon Y\to Y_1$ is defined by $D_2$, and both $f_1,f_2$ are birational since $q$ is birational (so $D_1,D_2$ are nef and big). 
Since $\rho(Y)=1$ and $D_2$ is nef and big, it follows that $D_2$ is ample and hence that $f_2$ is an isomorphism. 
\end{proof}

We are now ready for the proof of our main result.

\begin{theorem}\label{thm:not contain in max subgroups}
Let $n\geq 2$ and $m\geq0$ be a positive integers.
Let $X$ be a non-rational and rationally connected variety carrying a non-trivial conic bundle structure and admitting a fibration $c\colon X\to \IP^2$.
Set $\mathcal P_n=\IP_X(\mathcal O_X\oplus c^*\mathcal O_{\IP^2}(n))$.
Then the group $\Autz(\mathcal P_n\times\IP^m)$ is not contained in a maximal connected algebraic group of $\Bir(\mathcal P_n\times\IP^m)$. 
\end{theorem}
\begin{proof}
If $m=0$, the statement is Theorem~\ref{thm:main1new}, so let us assume that $m\geq1$. 
Notice that $\Autz(\mathcal P_n\times\IP^m)\simeq\Autz(\mathcal P_n)\times\Aut(\IP^m)$. 
Assume that $\Autz(\mathcal P_n\times\IP^m)$ is contained in a connected algebraic subgroup $H$ of $\Bir(\mathcal P_n\times\IP^m)$ acting rationally on $\mathcal P_n\times\IP^m$. 
By Weil regularisation theorem \cite{Wei55}, there is a birational map $\varphi\colon \mathcal P_n\times\IP^m\dashrightarrow V$ to a variety $V$ such that $H\subseteq \Autz(V)$. By Proposition~\ref{prop:bir maps between products},  we have $V\simeq W\times\IP^m$ and $\varphi=(\varphi_1,\varphi_2)$, where $\varphi_1\colon \cP_n\dashrightarrow W$ is birational and $\varphi_2$ is an isomorphism. 
Notice that $ \Autz(V)\simeq \Autz(W)\times\Aut(\IP^m)$. 
By Theorem~\ref{thm:main1new}, there exists a
variety $Y$ and an $\Autz(W)$-equivariant birational map $W\dashrightarrow Y$ such that $\Autz(W)\subsetneq \Autz(Y)=:G$. 
Then $H\subsetneq G\times\Aut(\IP^m)$. 
\end{proof}

\begin{proof}[Proof of Main Theorem]
The threefold $X$ from \cite[Example 2-6]{MB86} is irrational and has a fibration $X\to\IP^2$ and  $X\times\IP^3$ is rational \cite{BCTSSD85}. 
In \cite{SB02} it is shown that already $X\times\IP^2$ is rational. 
In particular, $\mathcal P_n\times\IP^m$ is rational for $n\geq0$ and $m\geq1$. 
The claim now follows from Theorem~\ref{thm:not contain in max subgroups} applied to the rational variety $Y:=\cP_n\times\IP^m$ for $n\geq2$ and $m\geq1$, which is of dimension $\dim Y=4+m\geq5$.
\end{proof}

\begin{remark}
We notice that if $X$ is any stably rational and non-rational variety of dimension 3,
and $k_0$ is the smalles positive integer such that $X\times\IP^{k_0}$ is rational, 
then by Theorem \ref{thm:not contain in max subgroups},
the group $\Autz(\mathcal P_n\times\IP^m)$ is not contained in a maximal connected algebraic group of $\Bir(\IP^{m+4})$ for any $m\in\IN$ such that $m+4\geq k_0+3$.
\end{remark}

\appendix

\section{}

We start with the definition of the $Quot$ functor.

\begin{definition}
Let $S$ be a noetherian scheme, $X$ a scheme of finite type over $S$ and $\mathcal E$ a coherent sheaf on $X$.
We define the functor $Quot_{\mathcal E/X/S}\colon (Sch/S)\to (SET)$
by setting on objects
$$Quot_{\mathcal E/X/S}(T)=\{\mathcal E_T\to \mathcal F\to 0\vert\;\mathcal F\text{ with proper support and flat over }T\}.$$
\end{definition}
If $X/S$ is projective, by a result of Grothendieck the functor $Quot_{\mathcal E/X/S}$ is representable by a locally projective scheme $\widetilde{\mathcal Q}$.
In this appendix we prove that, if there is a group $G$ acting on $X$ and $S$ making $X\to S$ $G$-equivariant and $\mathcal E$ is $G$-linearised,
then there is an action of $G$ on $\widetilde{\mathcal Q}$.

More precisely, we prove the following.

\begin{proposition}\label{prop:App}
Let $X/S$ be a projective scheme and $\mathcal E$ a coherent sheaf on $X$.
Let $G$ be a group.
Assume that there is an action of $G$  on $X$ and $S$ making $X\to S$ $G$-equivariant and $\mathcal E$ is $G$-linearised.
Let $\widetilde{\mathcal Q}$ be the $S$-scheme representing $Quot_{\mathcal E/X/S}$. Then
\begin{enumerate}
\item there is an action of $G$ on $\widetilde{\mathcal Q}$ making $\widetilde{\mathcal Q}\to S$ $G$-equivariant, and
\item for every $S$-scheme $U$ carrying a $G$-action making $U\to S$ $G$-equivariant, 
for every $G$-linearised quotient of $\mathcal E_U$, the induced morphism $U\to\widetilde{\mathcal Q} $ is $G$-equivariant.
\end{enumerate}
\end{proposition}

\begin{proof}
Let $g\in G$. Then $g$ induces an isomorphism of categories

$$
\begin{array}{rccc}
\Phi_g\colon&(Sch/S)&\rightarrow&(Sch/S)\\
&\gamma\colon T\to S&\mapsto&g\circ\gamma\colon T\to S
\end{array}
$$
\noindent We notice that $\Phi_{g_1\circ g_2}=\Phi_{g_1}\circ \Phi_{g_2}$.

\smallskip

\noindent For every $T$ object of $(Sch/S)$ we have
$$
Quot_{\mathcal E/X/S}(\Phi_g(T))=\{\mathcal{E}_{\Phi_g(T)}={(g^*\mathcal{E})_T\to} \mathcal F\to 0\vert\;\mathcal F\text{ flat over }T\},
$$
\noindent therefore there is an induced bijection
$$
\begin{array}{rccc}
F_g(T)\colon&Quot_{\mathcal E/X/S}(T)&\rightarrow&Quot_{\mathcal E/X/S}(\Phi_g(T))\\
&\mathcal E_T\to\mathcal F&\mapsto&\mathcal{E}_{\Phi_g(T)}\to\mathcal F.
\end{array}
$$

\noindent Moreover, if $G$ acts on $T$ and $T\to S$ is $G$-equivariant, then $F_g(T)(\mathcal F)=g^*\mathcal F$.

\smallskip

\noindent The bijection $F_g(T)$ induces a bijection
$$\widetilde F_g(T)\colon Hom_s(T, \widetilde{\mathcal Q})\to  Hom_s(\Phi_g(T), \widetilde{\mathcal Q}).$$
\noindent We set $\varphi_g=\widetilde F_g(\widetilde{\mathcal Q})(id_{\widetilde{\mathcal Q}})\in Hom_s(\Phi_g(\widetilde{\mathcal Q}), \widetilde{\mathcal Q})$.
Thus, $\varphi_g$ fits into a diagram

$$
\xymatrix{
\widetilde{\mathcal Q}\ar[r]^{\varphi_g}\ar[d]&\widetilde{\mathcal Q}\ar[d]\\
S\ar[r]_{g}&S
}
$$

\noindent Since $\Phi_{g_1\circ g_2}=\Phi_{g_1}\circ \Phi_{g_2}$, we have $\varphi_{g_1\circ g_2}=\varphi_{g_1}\circ \varphi_{g_2}$. This proves (1).

\medskip

The bijection $\widetilde F_g(T)$ is given by sending the diagram

$$
\xymatrix{
T\ar[rr]\ar[rd]&&\widetilde{\mathcal Q}\ar[ld]\\
&S&
}
$$

to

$$
\xymatrix{
T\ar[rr]\ar[rd]&&\widetilde{\mathcal Q}\ar[ld]\ar[rr]^{\varphi_g}&&\widetilde{\mathcal Q}\ar[lldd]\\
&S\ar[rd]_{g}&&&\\
&&S&&
}
$$

Moreover, if there is an action on $T$ making $T\to S$ equivariant,
then from the last diagram we deduce the diagram 
$$
\xymatrix{
T\ar[rrd]\ar[r]^{g^{-1}} &T\ar[rr]&&\widetilde{\mathcal Q}\ar[r]^{\varphi_g}&\widetilde{\mathcal Q}\ar[lld]\\
&&S&&
}
$$
Let $U, \mathcal F$ be as in (2).
Let $\alpha\colon U\to \widetilde{\mathcal Q}$ be the associate morphism.
Since $\mathcal F$ is $G$-linearized, the morphism induced by $F_g(U)(\mathcal F)$ is again $\alpha$.
Therefore, by our previous discussion, this is also $\varphi_g\circ\alpha\circ g^{-1}$,
proving the claim.
\end{proof}

\bibliographystyle{abbrv}
\bibliography{biblio}

\end{document}